\documentclass[12pt]{amsproc}

\usepackage[top=1in,bottom=1in,left=1.5in,right=1.5in]{geometry}
\usepackage{amsfonts} 
\usepackage{amsmath}
\usepackage{amssymb}
\usepackage{amsthm}
\usepackage{mathrsfs}
\usepackage{setspace}
\usepackage{subfigure}
\usepackage{url}
\usepackage{booktabs}
\usepackage{ifthen}
\usepackage{tikz}
\usetikzlibrary{calc}
\usepackage{enumerate}

\newcommand{\A}{\mathcal{A}}

\newcommand{\C}{\mathscr{C}}

\newtheorem{thm}{Theorem}[section]

\newtheorem{claim}[thm]{Claim}

\numberwithin{equation}{section}

\begin{document}

\title[Vizing's conjecture]{Vizing's conjecture for cographs}

\author{Elliot Krop}
\address{Elliot Krop (\tt elliotkrop@clayton.edu)}
\address{Department of Mathematics, Clayton State University}
\date{\today}

\begin {abstract}
We show that if $G$ is a cograph, that is $P_4$-free, then for any graph $H$, $\gamma(G\square H)\geq \gamma(G)\gamma(H)$. By the characterization of cographs as a finite sequence of unions and joins of $K_1$, this result easily follows from that of Bartsalkin and German. However, the techniques used are new and may be useful to prove other results.
\\[\baselineskip] 2010 Mathematics Subject
      Classification: 05C69
\\[\baselineskip]
\end {abstract}

\maketitle

 \section{Introduction}
Vizing's conjecture \cite{Vizing}, now open for fifty-three years, states that for any two graphs $G$ and $H$,
\begin{align}
\gamma(G \square H) \geq \gamma(G)\gamma(H)\label{V}
\end{align}
where $\gamma(G)$ is the domination number of $G$.

The survey \cite{BDGHHKR} discusses many results and approaches to the problem. For more recent partial results see \cite{ST}, \cite{PPS}, \cite{B}, \cite{CK}, \cite{K}, and \cite{K2}.

A predominant approach to the conjecture has been to show it true for some large class of graphs. For example, in their seminal result, Bartsalkin and German \cite{BG} showed the conjecture for decomposable graphs. More recently, Aharoni and Szab\'{o} \cite{AS} showed the conjecture for chordal graphs and Bre\v{s}ar \cite{B} gave a new proof of the conjecture for graphs $G$ with domination number $3$. 

We say that a bound is of \emph{Vizing-type} if $\gamma(G\square H)\geq c \gamma(G)\gamma(H)$ for some constant $c$, which may depend on $G$ or $H$. It is known \cite{ST} that all graphs satisfy the Vizing-type bound, 
\begin{align*}
\gamma(G \square H) \geq \frac{1}{2}\gamma(G)\gamma(H)+\frac{1}{2}\min\{\gamma(G),\gamma(H)\}.
\end{align*}

Restricting the graphs, but as a generalization of Bartsalkin and German's class of decomposable graphs, Contractor and Krop \cite{CK} showed
\begin{align*}
\gamma(G\square H)\geq \left(\gamma(G)-\sqrt{\gamma(G)}\right)\gamma(H)
\end{align*}

where $G$ belongs to $\A_1$, the class of graphs which are spanning subgraphs of domination critical graphs $G'$, so that $G$ and $G'$ have the same domination number and the clique partition number of $G'$ is one more than its domination number.

Krop \cite{K} showed that any claw-free graph $G$ satisfies the Vizing-type bound
\[\gamma(G\square H)\geq \frac{2}{3}\gamma(G)\gamma(H)\]

In this paper we show that the class of induced $P_4$-free graphs, or cographs, satisfies Vizing's conjecture.

\subsection{Notation}

All graphs $G(V,E)$ are finite, simple, connected, undirected graphs with vertex set $V$ and edge set $E$. We may refer to the vertex set and edge set of $G$ as $V(G)$ and $E(G)$, respectively.  For more on basic graph theoretic notation and definitions we refer to Diestel~\cite{Diest}. 
 
For any graph $G=(V,E)$, a subset $S\subseteq V$ \emph{dominates} $G$ if $N[S]=G$. The minimum cardinality of $S \subseteq V$, so that $S$ dominates $G$ is called the \emph{domination number} of $G$ and is denoted $\gamma(G)$. We call a dominating set that realizes the domination number a $\gamma$-set.

 The \emph{Cartesian product} of two graphs $G_1(V_1,E_1)$ and $G_2(V_2,E_2)$, denoted by $G_1 \square G_2$, is a graph with vertex set $V_1 \times V_2$ and edge set $E(G_1 \square G_2) = \{((u_1,v_1),(u_2,v_2)) : v_1=v_2 \mbox{ and } (u_1,u_2) \in E_1, \mbox{ or } u_1 = u_2 \mbox{ and } (v_1,v_2) \in E_2\}$.

A graph $G$ is a \emph{cograph} or \emph{$P_4$-free} if it contains no induced $P_4$ subgraph. 

Let $G$ be any graph and $S$ a subset of its vertices. Chellali et al. \cite{CHHM} defined $S$ to be a \emph{$[j,k]$-set} if for every vertex $v\in V-S, j\leq |N(v) \cap S|\leq k$. Clearly, a $[j,k]$-set is a dominating set. For $k\geq 1$, the \emph{$[1,k]$-domination number} of $G$, written $\gamma_{[1,k]}(G)$, is the minimum cardinality of a $[1,k]$-set in $G$. A $[1,k]$-set with cardinality $\gamma_{[1,k]}(G)$ is called a \emph{$\gamma_{[1,k]}(G)$-set}.

If $\Gamma=\{v_1,\dots, v_k\}$ is a minimum dominating set of $G$, then for any $i\in [k]$, define the set of \emph{private neighbors} for $v_i$, $P_i=\big\{v\in V(G)-\Gamma: N(v)\cap \Gamma = \{v_i\}\big\}$. For $S\subseteq [k]$, $|S|\geq 2$, we define the \emph{shared neighbors} of $\{v_i:i\in S\}$, $P_S=\big\{v\in V(G)-\Gamma: N(v)\cap \Gamma=\{v_i: i\in S\}\big\}$.

For any $S\subseteq [k]$, say $S=\{i_1,\dots, i_s\}$ where $s\geq 2$. We may write $P_S$ as $P_{\{i_1,\dots, i_s\}}$ or $P_{i_1,\dots, i_s}$ interchangeably.

For $i\in [k]$, let $Q_i=\{v_i\} \cup P_i$. We call $\mathcal{Q}=\{Q_1,\dots, Q_k\}$ the \emph{cells} of $G$. For any $I\subseteq [k]$, we write $Q_I=\bigcup_{i\in I}Q_i$ and call $\C(\cup_{i\in I}Q_i)=\bigcup_{i\in I}Q_i\cup\bigcup_{S\subseteq I}P_{S}$ the \emph{chamber} of $Q_I$. We may write this as $\C_{I}$.

For a vertex $h\in V(H)$, the \emph{$G$-fiber}, $G^h$, is the subgraph of $G\square H$ induced by $\{(g,h):g\in V(G)\}$. 

For a minimum dominating set $D$ of $G\square H$, we define $D^h=D\cap G^h$. Likewise, for any set $S\subseteq [k]$, $P_S^h=P_S \times \{h\}$, and for $i\in [k]$, $Q_i^h=Q_i\times \{h\}$. By $v_i^h$ we mean the vertex $(v_i,h)$. For any $I^h\subseteq [k]$, where $I^h$ represents the indices of some cells in $G$-fiber $G^h$, we write $\C_{I^h}$ to mean the chamber of $Q^h_{I^h}$, that is, the set $\bigcup_{i\in I^h}Q_i\cup\bigcup_{S\subseteq I^h}P^h_{S}$.

We may write $\{v_i:i\in I^h\}$ for $\{v^h_i:i\in I^h\}$ when it is clear from context that we are talking about vertices of $G\square H$ and not vertices of $G$.

For clarity, assume that our representation of $G\square H$ is with $G$ on the $x$-axis and $H$ on the $y$-axis.

Any vertex $v\in V(G)\times V(H)$ is \emph{vertically dominated} if $(\{v\}\times N_H[h])\cap D \neq \emptyset$ and \emph{vertically undominated}, otherwise. For $i\in [k]$ and $h\in V(H)$, we say that the cell $Q_i^h$ is \emph{vertically dominated} if $(Q_i\times N_H[h])\cap D\neq\emptyset$. A cell which is not vertically dominated is \emph{vertically undominated}. 

In our argument, we label vertices of a minimum dominating set $D$ of $G\square H$, by labels from $[k]$ so that for any $i\in [k]$, projecting the vertices labeled by $i$ onto $H$ produces a dominating set of $H$. We call a vertex $(x, h)\in D^h$ with the single label $i$, \emph{free}, if there exists another vertex $(y, h) \in D^h$, which is given the label $i$.

\section{Cographs}

\begin{thm}\label{cograph}
For any cograph $G$ and any graph $H$, $\gamma(G\square H)\geq \gamma(G)\gamma(H)$.
\end{thm}

\begin{proof}
Let $\Gamma=\{v_1,\dots, v_k\}$ be a minimum $[1,2]$ dominating set of $G$ and let $D$ be a minimum dominating set of $G\square H$. By the result of Chellali et al. \cite{CHHM} (Theorem $8$), $\gamma(G)=k$. Suppose $u\in V(G)-\{\Gamma\}$ is adjacent to two vertices of $\Gamma$, say $v_1$ and $v_2$. Notice that if neither $v_1$ nor $v_2$ have private neighbors with respect to $\Gamma$, then we could replace $v_1$ and $v_2$ by $u$ in $\Gamma$ and produce a smaller dominating set, which is a contradiction. Hence, at least one of $P_1$ or $P_2$ is nonempty.

\begin{claim}\label{external}
There exists a vertex in $P_1\cup P_2$ which is independent from both $u$ and $V(G)-\{Q_1\cup Q_2\}$.
\end{claim}
\begin{proof}
Case $1$: Suppose $P_1\neq \emptyset$ and $P_2=\emptyset$. Note that by the minimality of $\Gamma$ no vertex of $\Gamma-\{v_2\}$ can be adjacent to $v_2$. If $w_1\in P_1$, then by definition of private neighbors, no vertex of $\Gamma-\{v_1\}$ is adjacent to $w_1$. If $u$ is not adjacent to $w_1$, then we produce $P_4:w_1v_1uv_2$ which contradicts the definition of $G$. However, if $u$ is adjacent to every vertex of $P_1$, then we could replace $v_1$ and $v_2$ by $u$ in $\Gamma$ which would produce a smaller dominating set, which is impossible.


\begin{figure}[ht]
\begin{center}
\begin{tikzpicture}[scale=1]
\tikzstyle{vertex}=[circle, draw, inner sep=0pt, minimum size=6pt]
\tikzstyle{vert}=[circle,fill=black,inner sep=3pt]
\tikzstyle{overt}=[circle,fill=black!30, inner sep=3pt]
\tikzstyle{root}=[rectangle,fill=black,inner sep=3pt]

  \node[root, label=below:\tiny{$v_1$}] (v1) at (1,1) {};
  \node[root, label=below:\tiny{$v_2$}] (v2) at (2,1) {};

  \node[vertex, label=above:\tiny{$w_1$}](w1) at (0.5,2){};
  \node[vertex, label=above:\tiny{$u$}](u) at (1.5,2){};

 \draw[color=black] 
(w1)--(v1)--(u)--(v2)
;

\end{tikzpicture}
\caption{}
\end{center}

\end {figure}

Case $2$: Suppose $P_1,P_2\neq \emptyset$. By minimality of $\Gamma$, some vertex of $P_1\cup P_2$ is not adjacent to $u$. Suppose such a vertex is $w_2\in P_2$. We may assume $v_1$ is adjacent to $v_2$, else we would produce $P_4: w_2v_2uv_1$. For any vertex $w_1\in P_1$, we may assume $w_1$ is adjacent to $w_2$, else we would produce $P_4: w_1v_1v_2w_2$. Notice $u$ is adjacent to $w_1$ to avoid $P_4: w_2w_1v_1u$. Suppose $u'\in V(G)-\{Q_1\cup Q_2\}$ is adjacent to $w_2$ and suppose without loss of generality that $u'$ is adjacent to $v_3$. Thus, we are left with the situation illustrated in Figure $2$ where the drawn edges have been shown to exist.

\begin{figure}[ht]
\begin{center}
\begin{tikzpicture}[scale=1]
\tikzstyle{vertex}=[circle, draw, inner sep=0pt, minimum size=6pt]
\tikzstyle{vert}=[circle,fill=black,inner sep=3pt]
\tikzstyle{overt}=[circle,fill=black!30, inner sep=3pt]
\tikzstyle{root}=[rectangle,fill=black,inner sep=3pt]

  \node[root, label=below:\tiny{$v_1$}] (v1) at (1,1) {};
  \node[root, label=below:\tiny{$v_2$}] (v2) at (2,1) {};
  \node[root, label=below:\tiny{$v_3$}] (v3) at (3,1) {};

  \node[vertex, label=above:\tiny{$w_1$}](w1) at (0.5,2){};
  \node[vertex, label=above:\tiny{$u$}](u) at (1.5,2){};
  \node[vertex, label=above:\tiny{$u'$}](u') at (3,2){};
  \node[vertex, label=above:\tiny{$w_2$}](w2) at (2.5,2){};

 \draw[color=black] 
   (v3)--(u')--(w2)--(v2)  (w1)--(v1)--(u)--(v2)--(v1) (u)--(w1)

   (w1).. controls (1.5,2.25) ..(w2);

\end{tikzpicture}
\caption{}
\end{center}

\end {figure}

Since $u'$ may adjacent to at most two vertices of $\Gamma$ we argue that $u'$ is adjacent to either $v_1$ or $v_2$, since otherwise we have $P_4: u'w_2v_2v_1$. 

Subcase $(i)$: If $u'$ is adjacent to $v_2$, then $u'$ is also adjacent to $w_1$ to avoid $P_4: v_3u'w_2w_1$. Furthermore, if $v_3$ is not adjacent to $v_1$ or $v_1$, then we produce $P_4: v_3u'v_2v_1$. If $v_3$ is adjacent to $v_1$, then we produce $P_4: w_2v_2v_1v_3$. Thus, $v_2$ is adjacent to $v_3$. However, now we have $P_4: v_3v_2uw_1$ which is impossible.

Subcase $(ii)$: If $u'$ is adjacent to $v_1$, then $u'$ is also adjacent to $w_1$ to avoid $P_4: v_3u'w_2w_1$. Furthermore, $v_3$ is adjacent to $v_2$, else we have $P_4: v_3u'w_2v_2$. However, this forces $P_4: v_3v_2uw_1$ which is impossible.

\end{proof}

For any $h\in V(H)$, suppose the fiber $G^h$ contains $\ell_h(=\ell)$ vertically undominated cells $U^h=\big\{Q_{i_1}^h,\dots, Q_{i_{\ell}}^h\big\}$ for $0\leq \ell \leq k$. We set $I^h=\{i_1,\dots, i_{\ell}\}$. Notice that for $j_1,j_2\in [k]-I^h$, no vertex of $P^h_{j_1,j_2}$ may dominate any of $v_{i_1},\dots,v_{i_{\ell}}$. Thus, $\{v_i:i\in I^h\}$ must be dominated horizontally in $G^h$ either by shared neighbors of $\{v_i:i\in I^h\}$ or by vertices of $\{v_i:1\leq i \leq k, i\notin I^h\}$. Furthermore, the private neighbors $\{P_i^h:i\in I^h\}$ must be dominated horizontally in $G^h$ either by shared neighbors of $\{v_i:i\in I^h\}$ or by vertices of $\{P_i^h:1\leq i \leq k, i\notin I^h\}$.

\medskip

We label the vertices of $D$ by the following \emph{Provisional Labeling}. If a vertex of $D^h$ for any $h\in H$, is in $Q_i^h$ for $1\leq i \leq k$, then we label that vertex by $i$. If $v$ is a shared neighbor of some subset of $\{v_i:i\in I^h\}$, then it is a member of $P^h_{i,j}$ for some $i,j\in I^h$, and we label $v$ by the pair of labels $(i,j)$. If $v$ is a member of $P^h_{i,j}$ for $i\in I^h$ and $j\in [k]-I^h$, then we label $v$ by $i$. If $v$ is a member of $P^h_{i,j}$ for $i,j \in [k]-I^h$, then we label $v$ by either $i$ or $j$ arbitrarily.

After the labels are placed, all vertices of $D$ have a single label or a pair of labels.

\medskip

Next, we apply a relabeling to some of the vertices of $D$ which we call the \emph{First Refinement}. For a fixed $h\in H$, suppose $v$ is some shared neighbor of two vertices of $\{v_i: i\in I^h\}$ in the chamber of $Q_{I^h}^h$, which is vertically dominated, say by $y\in D^{h'}$ for some $h'\in H,\,h\neq h'$. In other words, we suppose $v\in P^h_{j_1,j_2}$ for some $j_1,j_2\in I^h$ which implies that $y\in P^{h'}_{j_1,j_2}$. 

The vertex $y$ may be labeled by one or two labels, regardless of whether the First Refinement had been performed on $D^{h'}$.

Suppose that $y$ is labeled by one label, say $j_1$. If $D^h$ contains a vertex $x\in P_{j_1,j_2}^h$, then we remove the pair of labels $(j_1,j_2)$ from $x$ and relabel $x$ by $j_2$. 

Suppose $y$ is labeled by the pair of labels, $(j_1,j_2)$. If $D^h$ contains a vertex $x\in P_{j_1,j_2}^h$, then we remove the pair of labels $(j_1,j_2)$ from $x$ and then relabel $x$ arbitrarily by one of the single labels $j_1$ or $j_2$.

After the labeling, a vertex $v$ of $D$ may have a pair of labels $(i,j)$ if $v\in P^h_{i,j}$ and for any $h'\in N_H(h)$, $D^{h'}\cap P^{h'}_{i,j}=\emptyset$.

\medskip

Finally, we relabel some of the vertices of $D$ by the \emph{Second Refinement}. For every $h\in H$, if $D^h$ contains vertices $x$ and $y$ with pairs of labels $(j_1,j_2), (j_2,j_3)$ respectively, for some integers $j_1,j_2,$ and $j_3$, then we relabel $y$ by the label $j_3$. If $x$ and $y$ are labeled $j_1$ and $(j_1,j_2)$ respectively, for some integers $j_1,j_2$, we relabel $y$ by $j_2$. We apply this relabeling to pairs of vertices of $D^h$, sequentially, in any order.

\medskip

\begin{claim}\label{single}
After the Second Refinement every label on a vertex of $D$ is a single label.
\end{claim}

\begin{proof}
For any $h\in V(H)$, suppose $v\in D^h$ has a pair of labels $(i,j)$. The Provisional Labeling prescribes that $i,j\in I^h$ which means that $Q_i$ and $Q_j$ are vertically undominated cells.  If there exists $w \in D^h \cap P_{j,m}$ for any $1\leq m \leq k$, or $x\in D^h\cap P_j$, then $v$ would have a single label after the Second Refinement which is not the case. By Claim \ref{external}, some vertex $x$ in $P_i^h\cup P_j^h$ is independent from $v$ and independent from $V(G)-\{Q_i\cup Q_j\}$. However, this means that $x$ is undominated, which contradicts the fact that $D$ is a dominating set.
\end{proof}

Suppose that for some $h\in V(H)$, $G^h$ contains a cell, $Q_i^h$, which is vertically undominated and the vertices of $D^h$ dominating $Q_i^h$ are not labeled $i$. in this case, $v_i^h$ can only be dominated by other members of $\{v_j^h:j\in [k], j\neq i\}$, so suppose for some $j_1\neq i, j_1\in [k]$, there exists $v_{j_1}^h\in D^h$ so that $v_i$ is adjacent to $v_{j_1}$. To avoid a contradiction to the minimality of $\Gamma$, we see that $P_i\neq \emptyset$ and say $u\in P_i^h$. Notice that if $u$ is dominated by some $u'\in P_{j_2}\cap D^h$ for some $j_2\neq i,j_1$, then we produce $P_4: v_{j_1}v_iuu'$ in $G^h$ and thus in $G$. Furthermore, if $v\in P_{j_1}\cap D^h$ dominates $u$, then $v$ is a free vertex labeled $j_1$ and we may relabel $v$ by $i$ without changing the vertically dominated status of cells $Q_{j_1}^{h'}$ for any $h'\in V(H)$. Finally, suppose $u$ is dominated by some shared neighbor $w\in P_{j_1,j_2}\cap D^h$. Notice if $x\in P_{j_1}$, then we produce $P_4: xv_{j_1}v_iv_{j_2}$ and if $y\in P_{j_2}$, then we produce $P_4: yv_{j_2}v_iv_{j_1}$ which cannot occur. Thus, $P_{j_1}=P_{j_2}=\emptyset$. If for every $w'\in P_{j_1,j_2}$, $w$ is adjacent to $w'$, we have a contradiction to the minimality of $\Gamma$, since now we can replace $v_{j_1}$ and $v_{j_2}$ by the projection of $w$ onto $V(G)$ and form a smaller dominating set of $G$. We are left to suppose there exists $w'\in P_{j_1,j_2}$ so that $w'$ is not adjacent to $w$. To avoid $P_4: w'v_{j_1}wu$, $w'$ must be adjacent to $u$. Furthermore, this same property is true for any vertex $v\in P_{j_1,j_2}$ not adjacent to $w'$ or $w$, namely, $v$ must be adjacent to $u$. To avoid $P_4: v_{j_2}w'v_{j_1}v_i$ we must also have $v_{j_2}$ adjacent to $v_i$. At this point, notice that \{$v_{j_1},v_{j_2}, P_{j_1,j_2}$\} is dominated by $v_i$ and $u$, which is a contradiction to the minimality of $\Gamma$, since now we can replace $v_{j_1}$ and $v_{j_2}$ in $\Gamma$ by the projection of $u$ onto $V(G)$ and produce a smaller dominating set of $G$.

Notice that for any $i\in [k]$, projecting all vertices with a label $i$ to $H$ produces a dominating set of $H$. Summing over all $i$, we count at least $\gamma(G)\gamma(H)$ vertices in $D$.

\end{proof}

 \bibliographystyle{plain}
 
 \end{document}